\documentclass[reqno]{amsart}
\usepackage{hyperref}
\begin{document}
\title[Oscillation of  DDEs via the hyper4 convergence]{Oscillation of  delay  differential equations via the hyper4 convergence}
\author[G. L. Karakostas]
{George L. Karakostas}
\address{George L. Karakostas \newline
 Department of Mathematics, University of Ioannina,
 451 10 Ioannina, Greece}
\email{gkarako@uoi.gr}
\date{}
 \begin{abstract} A sharp condition is  provided to guarantee that the (nontrivial) solutions of a DDE of the form $\dot{x}(t)+F(t,x)=0$ $t\geq 0,$ (where $F(t,\cdot)$ is an odd-like causal operator) either oscillate, or converge monotonically to zero. The method used is based on the convergence of the sequence of hyper4-iterations to the Lambert's function.
 \end{abstract}
 \maketitle
\numberwithin{equation}{section}
\newtheorem{theorem}{Theorem}[section]
\newtheorem{lemma}[theorem]{Lemma}
\newtheorem{definition}[theorem]{Definition}
\newtheorem{example}[theorem]{Example}
\newtheorem{property}[theorem]{Property}
\newtheorem{remark}[theorem]{Remark}
\newtheorem{application}[theorem]{Application}
\section{Introduction}
Theorem 2.3.1 in \cite{GL} states as follows: {\em Consider the linear non-autonomous delay differential equation
\begin{equation}\label{e1}\dot{x}(t)+P(t)x(t-\tau)=0,\end{equation} 
$t\geq t_0.$ Assume that $P\in C([t_0,\infty),\mathbb{R}^+),$ $\tau>0$ and $$\liminf_{t\to+\infty}\int_{t-\tau}^tP(s)ds>\frac{1}{e}.$$ Then every solution oscillates.}\par   
In this short notice we extend the previous  oscillation result to more general delay differential equations by using a method which lies on the range of the convergence of the sequence of tetrations to the Lambert's function.  \par Let  $\sigma<0$ and  $T:C(I, \mathbb{R})\to C(I, \mathbb{R})$  be an operator, where $I:=(\sigma, +\infty)$. We say  that $T$ has an {\em amnesia on $I$}, if for each $t\geq 0$, there is a subset $H(t)\subseteq (\sigma,t)$ with the property that, whenever $u, v\in C(I, \mathbb{R})$ are such that $u(r)=v(r),$ for all $r\in H(t)$,  then $(Tu)(t)=(Tv)(t).$  Let $\sigma(t):=\inf H(t)$ and $\tau(t):=\sup H(t).$ The interval $I(t):=[\tau(t), t]$ is the {\em amnesia interval} at $t$ of the operator. Differential equations with response an operator with amnesia have specific asymptotic properties. Some of them are exhibited elsewhere, see \cite{KS1, KS2}.
In this paper we give sufficient conditions for oscillation of the solutions of a differential equation of the form \begin{equation}\label{e2}\dot{u}(t)+(Tu)(t)=0, \enskip t\geq 0,\end{equation}
where, $T$ is an operator with amnesia on $I$. \par
  Notice that a great number of authors studied  the linear differential equation $\dot{x}(t)+p(t)x(\tau(t))=0,$ $t\geq t_0$, under the condition $\liminf_{t\to+\infty}\int_{t-\tau}^tP(s)ds\leq 1/e$ see, e.g.\cite{D1},  \cite{GL}, \cite{EKZ} - \cite{S}. The discrete version of equation \eqref{e2} is studied  in \cite{G1}.
\par 
Before presenting our results we shall repeat the basic meanings of the so called tetration, or hyper4, which is a hyperoperation of the form $x^{x^{x^{x^{\cdots}}}}$ ($n$ exponents) which, according to the Hooshmand's notation \cite{HOO}, is denoted as $x^{\underline{n}}$. It is known that  the limit $y$ of $x^{\underline{n}}$, as $n$ goes to infinity, exists for the values of $x$ in the interval $e^{-e}\leq x \leq e^{1/e}$, which approximately is the  interval [0.065988, 1.4446678], a result shown by  Euler, \cite{EU}. Thus, $x = y^{1/y}$. The limit defining the infinite exponential of $x$ does not exist when $x > e^{1/e}$ because the maximum of $y^{1/y}$  is 
$e^{1/e}$. The limit also fails to exist when $0 < x < e^{-e}$.
It is obvious that the limit, if it exists, is a positive real solution of the equation $y = x^y$. Moreover in this case we observe that the relation
$y=x^y=e^{y\ln{x}}$ holds if and only if 
$-\ln{x}=(-y\ln{x})e^{-y\ln{x}}.$ This implies that if  $W$ represents the well known Lambert's function, i.e.,  the function $W:=W(x)$ defined by $x=We^W,$ then we have $W(-\ln{x})=-y\ln{x}$ and so 
$$\frac{W(-\ln{x})}{-\ln{x}}=y.$$  
\par 
In section 2 of the paper we state and prove the main oscillation result and in the last section 3
some simple applications illustrate the result. 
\section{The main results}  
Consider the differential equation \eqref{e2} and assume that  the response operator satisfies the following basic condition:\par (C)  There is a positive Lebesgue integral function $b(t)$ such that 
\begin{equation}\label{e4}(Tx)(t)\geq b(t)\inf_{r\in I(t)}x(r),\enskip\text{whenever}\enskip\inf_{r\in I(t)}x(r)>0\end{equation} and 
 \begin{equation}\label{e5}(Tx)(t)\leq b(t)\sup_{r\in I(t)}x(r),\enskip\text{whenever}\enskip\sup_{r\in I(t)}x(r)<0,\end{equation} where, recall that, $I(t)$ is the amnesia interval at $t$.
  It is clear that if $T$ is a positive linear operator, then it satisfies condition (C). \par These two mathematical relations can be written as 
  $$\min\{(Tx)(t), -(T(-x))(t)\}\geq b(t)\inf_{r\in I(t)}x(r), \enskip\text{whenever}\enskip\min_{r\in I_n}x(r)>0$$
and, if $(T\cdot)$ is odd, then conditions \eqref{e4} and \eqref{e5}  are equivalent. \par 
  Our result if this work is given in the following theorem, where by the property $(\bf{P})$ we shall mean that a function oscillates or it tends monotonically to zero. :
 
 \begin{theorem}\label{t2} Assume that condition (C)  holds, as well as that \begin{equation}\label{e6} \liminf_{t+\to +\infty} \sigma(t)=+\infty.\end{equation} 
If it holds \begin{equation}\label{e8}w:=\liminf_{t\to+\infty}\int_{\tau(t)}^tb(s)ds> \frac{1}{e},\end{equation}  then property {\rm(\textbf{P})}  keeps in force. \end{theorem}\par 
\begin{proof} Assume that $x$ is a (nontrivial) non-oscillating solution of equation \eqref{e2}. If $x$ is eventually negative, the second relation in (C) is satisfied for all  large $t$. Then the function $y:=-x$ is eventually positive and it satisfies the differential equation $\dot{y}(t)+(T^*y)(t)=0,$ where the new function $(T^*u):=-T(-u)$, obviously,  satisfies the first condition in (C) for all large $t.$ Thus,  we can assume that $x$ is eventually positive, i.e. there is some $t_0$ such that $x(r)>0,$ for all $r\geq t_0.$ Due to \eqref{e6}, there is some $\bar{t}\geq t_0$ such  that $t_0\leq \sigma(t),$  for all $t\geq \bar{t}.$ Then, for all $t\geq \bar{t}$ and $r\in I(t)$, we have $x(r)>0$, and therefore 
\begin{equation}\label{e9}\dot{x}(t)=-(Tx)(t)< 0.\end{equation} 
Thus $x$ is an eventually monotonically decreasing function and so $$x(\sigma(t))\geq x(r)\geq x(\tau(t))>x(t), \enskip t\geq \bar{t}\enskip\text{and}\enskip r\in I(t).$$\par 
Due to  condition (C) we get
\begin{equation}\label{e10}\dot{x}(t)=-(Tx)(t)\leq-b(t)x(\tau(t)),\enskip t\geq \bar{t}.\end{equation}
  \par Let us define $$\zeta:=\liminf_{t\to\infty}\frac{x(\tau(t))}{x(t)},$$ which, in case the solution does not converge to zero, is bounded above. Obviously, $\zeta\geq 1.$ Let $0<\epsilon<1$. Then for some $\hat{t}>\bar{t}$ and all $t\geq \hat{t}$ it holds
$$\frac{x(\tau(t))}{x(t)}>\zeta-\epsilon.$$ Hence for all $t\geq \hat{t}$ we have
$$\frac{\dot{x}(t)}{x(t)}=-\frac{1}{x(t)}(Tx)(t)\leq-\frac{x(\tau(t))}{x(t)}b(t)\leq -(\zeta-\epsilon)b(t),$$ and therefore
$$\ln\frac{x(t)}{x(\tau(t))}\leq-(\zeta-\epsilon)\int_{\tau(t)}^tb(s)ds.$$ Also, there is a sequence $(t_n)$ converging to $+\infty$ such that
$$\zeta+\epsilon>\frac{x(\tau(t_n))}{x(t_n)}.$$  
We can assume that $t_n\geq \hat{t},$ for all $n=1, 2, \cdots.$  Therefore we have 
$$\zeta+\epsilon\geq  e^{(\zeta-\epsilon)\int_{\tau(t_n)}^{t_n}b(s)ds}, \enskip n=1, 2, \cdots.$$ 
From the previous results it follows that
$$\zeta+\epsilon\geq \liminf_{t\to+\infty}e^{(\zeta-\epsilon)\int_{\tau(t)}^{t}b(s)ds}\hskip .9in $$
$$\hskip.3in =e^{(\zeta-\epsilon)\liminf_{t\to+\infty}\int_{\tau(t)}^{t}b(s)ds}=e^{(\zeta-\epsilon)w}$$ and taking the limits of both sides when $\epsilon$ tends to zero, we get
\begin{equation}\label{e11}\zeta\geq e^{\zeta w}.\end{equation} 
\par
Since it holds $\zeta\geq 1$, from \eqref{e11} we obtain $\zeta\geq e^w=:a.$ Then, again, we have 
$\zeta\geq e^{\zeta w}=(e^w)^{\zeta}=a^{\zeta}\geq a^a$
and inductively it follows
$$\zeta\geq a^{\underline{n}},$$ for all $n=1, 2, \cdots.$
We observe that $1<e^w=a$  and by induction $a^{\underline{n-1}}<a^{\underline{n}},$ namely the sequence  of tetrations $(a^{\underline{n}})$ is increasing and it has an upper bound, the real number $\zeta.$  Thus it must converge.  But, as we said in the introduction,  this convergence  may occur only in case
$${e^{-e}}\leq a\leq e^{e^{-1}},$$ which implies that $w\leq e^{-1}.$ Obviously, this is impossible, due to \eqref{e8}.\par 
The proof is complete.
  \end{proof}

 \section{Applications}
{\bf Application 1}:    Consider the difference equation \par
$$\dot{x}(t)+\frac{1}{qt}x(t-6)+\frac{t-1}{qt}x(t-8)=0,\enskip t\geq 6,$$ where $q>0.$
Here  we have $\tau(t)=t-6$ and $b(t)=\frac{1}{q}$. Clearly, if $q>6e,$ then condition  \eqref{e8} is true. Thus for these values of $q$ all solutions have property (\rm{\textbf P}).\par
{\bf Application 2}:   Let  $a_2, a_3>0,$ $a:=\min\{a_2,a_3\}$ and let $a_1\in\mathbb{R}\setminus\{0\}$ such that $ae^{1+a_1}(e^{a_1}-1)/a_1>1.$ Consider  the differential equation  $$\dot{x}(t)+\int_1^2e^{\max\{a_1s,x^2(t-a_2s)\}}x(t-a_3s)ds=0,\enskip t\geq 0.$$ Here we have $\tau(t):=t-a$, where and $b(t)=e^{a_1}(e^{a_1}-1)/a_1.$ Now, condition  \eqref{e8} is satisfied by the choice of the parameters $a_1, a_2, a_3.$ Thus all solutions have property (\rm{\textbf P}). \par {\bf Application 3}: Consider the difference equation  $$\dot{x}(t)+\int_0^1[as^m+bs^2\sin^l(x^3(t-s-5))]x(t-s-1)ds=0,\enskip t\geq 0,$$ where $a, b, m$ are positive numbers such that $(a-mb)e>m$ when the positive integer $l$ is odd, and $ae>m,$ whenever $l$ is positive even integer. In this case we have $\tau(t)=t-1$ and $b(t)=\frac{a}{m}-b,$ when $l$ is odd and $b(t)=\frac{a}{m},$ when $l$ is even. These conditions guarantee the fact that all solutions have the property ({\bf P}). 
\par 
 

\begin{thebibliography}{XX}
  \bibitem{D1} Julio G. Dix, Improved oscillation criteria for first-order delay differential equations with variable delay, {\em Electron. J. Diff, equations} 32 (2021), 1-12.
    \bibitem{EKZ} L.H. Erbe, Qingkai Kong, B.G. Zhang, {\em Oscillation Theory for Functional Differential Equations}, Marcel Dekker, New York, 1995. 
  \bibitem{EZ} L.H. Erbe, B.G. Zhang, Oscillation of first order linear differential equations with deviating arguments, {\em Differ. Integral Equ.} 1 (1988) 305–314. 
  \bibitem{EU} L. Euler,  De serie Lambertina Plurimisque eius insignibus proprietatibus. {\em  Acta Acad. Scient. Petropol.} 2, 29–51, 1783. Reprinted in Euler, L. Opera Omnia, Series Prima, Vol. 6: Commentationes Algebraicae. Leipzig, Germany: Teubner, pp. 350–369, 1921. 
  \bibitem{FK} N. Fukagai, T. Kusano, Oscillation theory of first order functional differential equations with deviating arguments, {\em Ann. Mat. Pura Appl.} 136 (1984) 95-117.
 \bibitem{GL} I. Gy\"ori and G. Ladas, {\em Oscillation Theory of Delay Differential Equations With Applications}, Clarendon Press, Oxforf, 1991.
 \bibitem{HOO} Hooshmand, M. H. (2006). Ultra power and ultra exponential functions, {\em Integral Transforms and Special Functions}. 17 (8): 549–558.
 \bibitem{WI} https://en.wikipedia.org/wiki/Tetration (download: Dec. 07, 2021)
   \bibitem{JS} J. Jaro\^s, I.P. Stavroulakis, Oscillation tests for delay equations, {\em Rocky Mountain J. Math.} 29 (1999) 139–145. 
 \bibitem{G1} G. L. Karakostas, Oscillation of  difference  causal operator equations, {\em Journal of Computational Mathematica }, 2456-8686, 6(1), 2022: 091-106. https://doi.org/10.26524/cm123.
 \bibitem{KS1} G. Karakostas and Y. Sficas, Uniform Asymptotic Stability of Delay Differential Equations with Amnesia, {\em J. Math. Anal. Appl.} 194 (1995),  437-458.
   \bibitem{KS2} G. Karakostas and Y. Sficas, Stability of Delay Differential Equations with Amnesia and a Sign-Guiding Factor {\em J. Math. Anal. Appl.} 220 (1998),  204-223.
 \bibitem{KSS} M. Kon, Y.G. Sficas, I.P. Stavroulakis, Oscillation criteria for delay equations, {\em Proc. Amer. Math. Soc.} 128 (2000) 2989–2997. 
  \bibitem{KC} R.G. Koplatadze, T.A. Chanturija, On the oscillatory and monotonic solutions of first order differential equations with deviating arguments, 
{\em Differentsial ’nye Uravneniya} 18 (1982) 1463–1465.
 \bibitem{K} M.K. Kwong, Oscillation of first order delay equations, {\em J. Math. Anal. Appl.} 156 (1991) 374. 286.
\bibitem{L1} G. Ladas, Sharp conditions for oscillations caused by delay, {\em Applicable Anal.} 9 (1979) 93–98.
 \bibitem{LLZ} G.S. Ladde, V. Lakshmikantham, B.G. Zhang, {\em Oscillation Theory of Differential Equations with Deviating Arguments,} Marcel Dekker, New York, 1987. 
  \bibitem{L} B. Li, Oscillations of first order delay differential equations, {\em Proc. Amer. Math. Soc.} 124 (1996) 3729–3737. 
   \bibitem{S} I.P. Stavroulakis, Oscillation criteria for delay and difference equations with 
non-monotone arguments, {\em Applied Math.  Comput.} 226 (2014) 661–672
  \bibitem{SS}Y.G. Sficas, I.P. Stavroulakis, Oscillation criteria for first-order delay equations, {\em Bull. London Math. Soc.} 35 (2003) 239–246. 
   \bibitem{W} Z.C. Wang, I.P. Stavroulakis, X.Z. Qian, A Survey on the oscillation of solutions of first order linear differential equations with deviating arguments, {\em Appl. 
Math. E-Notes} 2 (2002) 171–191.
   \bibitem{YWZQ}  J.S. Yu, Z.C. Wang, B.G. Zhang, X.Z. Qian, Oscillations of differential equations with deviating arguments, {\em PanAmerican Math. J.} 2 (1992) 59–78. 
  \bibitem{Z} Y. Zhou, Y.H. Yu, On the oscillation of solutions of first order differential equations with deviating arguments, {\em Acta Math. Appl. Sinica} 15 (3) (1999), 288–302. 


 \end{thebibliography}
\end{document}